\documentclass{article}

\usepackage{arxiv}

\usepackage[utf8]{inputenc} %
\usepackage[T1]{fontenc}    %
\usepackage{xcolor}
\usepackage{hyperref}       %
\usepackage{url}            %
\usepackage{booktabs}       %
\usepackage{amsfonts}       %
\usepackage{nicefrac}       %
\usepackage{microtype}      %
\usepackage{lipsum}
\usepackage{graphicx}

\usepackage{amsmath}
\DeclareMathOperator*{\argmax}{\text{argmax}}

\usepackage{algorithm}
\usepackage{algpseudocode}
\usepackage{subcaption}
\usepackage{multirow}
\usepackage[square,numbers]{natbib}
\usepackage{doi}

\usepackage[normalem]{ulem}

\usepackage{caption} 
\title{Automated solver selection for simulation of multiphysics processes in porous media}

\begin{document}

\author{ \href{https://orcid.org/0009-0006-7095-3044}{\includegraphics[scale=0.06]{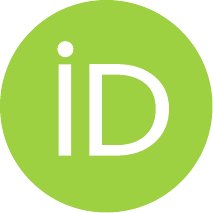}\hspace{1mm}Yury Zabegaev}\thanks{Corresponding author. Email: \href{mailto:yury.zabegaev@uib.no}{\texttt{yury.zabegaev@uib.no}}} \\
	Department of Mathematics\\
	University of Bergen\\
	Bergen, Norway \\
        \And
	\href{https://orcid.org/0000-0002-0333-9507}{\includegraphics[scale=0.06]{orcid.pdf}\hspace{1mm}Eirik Keilegavlen} \\
	Department of Mathematics\\
        University of Bergen\\
	Bergen, Norway \\
        \And
	Einar Iversen \\
	Department of Earth Science\\
        University of Bergen\\
	Bergen, Norway \\
        \And
	\href{https://orcid.org/0000-0002-0212-7959}{\includegraphics[scale=0.06]{orcid.pdf}\hspace{1mm}Inga Berre} \\
	Department of Mathematics\\
        University of Bergen\\
	Bergen, Norway \\
}

\date{}

\maketitle

\begin{abstract}
Porous media processes involve various physical phenomena such as mechanical deformation, transport, and fluid flow. Accurate simulations must capture the strong couplings between these phenomena. Choosing an efficient solver for the multiphysics problem usually entails the decoupling into subproblems related to separate physical phenomena. Then, the suitable solvers for each subproblem and the iteration scheme must be chosen. 
The wide range of options for the solver components makes finding the optimum difficult and time-consuming;
moreover, solvers come with numerical parameters that need to be optimized. 
As a further complication, the solver performance may depend on the physical regime of the simulation model, which may vary with time.
Switching a solver with respect to the dominant process can be beneficial, but the threshold of when to switch solver is unclear and complicated to analyze.
We address this challenge by developing a machine learning framework that automatically searches for the optimal solver for a given multiphysics simulation setup, based on statistical data from previously solved problems. 
For a series of problems, exemplified by successive
time steps in a time-dependent simulation, the framework updates and improves its decision model online during the simulation. We describe the solver selection algorithm, present examples of how the solver selector tunes the solver during the simulation, and show how it outperforms preselected state-of-the-art solvers for test problem setups. 
The examples are based on simulations of poromechanics and simulations of flow and transport. For the quasi-static linear Biot model, we demonstrate automated tuning of numerical solver parameters by showing how the L-parameter of the so-called Fixed-Stress preconditioner can be optimized.
Motivated by a test example where the main heat transfer mechanism changes between convection and diffusion, we also discuss how the solver selector can dynamically switch solvers when the dominant physical phenomenon changes with time.
\end{abstract}

\section{Introduction}

Processes in porous media involve strong multiphysics couplings, such as between (subsets of) mechanical deformation, fluid flow, and transport.
For such multiphysics simulations, a large part of the computational cost is associated with solving systems of linear equations.
Popular linear solvers for multiphysics problems apply preconditioners that split the linear system into subproblems related to separate physical phenomena. 
This allows for tailoring of the linear solver, as illustrated below for our model problems of poromechanics and coupled flow and transport, both of which are highly relevant for porous media applications.

Solvers for poromechanics have been an active research field for decades \cite{settari1998coupled}.
While direct solvers can be applied for small problems, constraints in CPU and memory necessitate the application of iterative solvers for application-relevant problem sizes. 
A common strategy is to apply block preconditioners; however, a naive block-Jacobi style splitting into subproblems is unstable and stabilization is needed to ensure efficient preconditioning; see e.g. \cite{White2016,both2017robust}.
This has led to the development of several splitting strategies \cite{kim2011stabilitydrained,kim2011stabilityfixed}, and further to the need to tune parameters within the splitting scheme \cite{Mikelic_Wheeler_fixed_stress,Storvik2018}.
Finally, the block preconditioner could be accelerated by an outer Krylov method.

Similar to poromechanics, coupled flow and transport in porous media involve interaction between processes with different mathematical structures \cite{trangenstein1989mathematical}:
The pressure field that governs the flow is parabolic. In most cases, it is almost elliptical.
Transport can involve components and/or energy and involve one or multiple fluid phases. It is commonly represented by an advection-diffusion equation, although one of the transport terms may be ignored depending on the application.
In subsurface engineering, transport is usually dominated by forced advection, and thus the problem is essentially hyperbolic.
Accordingly, the dominant preconditioning strategy, known as Constrained Pressure Residual (CPR) \cite{Wallis1983,Wallis1985}, involves two stages:
First, a pressure equation is derived by a decoupling technique \cite{Lacroix2003,stuben2007algebraic} and solved with an approach that handles the long-range couplings due to the elliptic components of the solution, commonly by a multigrid method, e.g. \cite{cao2005parallel,gries2014preconditioning}.
Second, a local correction is applied, often in the form of an incomplete lower-upper (ILU) factorization.
However, while CPR has been highly successful for transport problems that are essentially hyperbolic, it is less effective when diffusion plays an important role, as will be the case for many thermal problems \cite{li2015parallel}.
In this case, an alternative approach that acknowledges the parabolic nature of the transport problem has been shown to be more efficient; see \cite{Roy2019} and also \cite{cremon2020multi,cremon2023constrained}.
Independent of which strategy is applied, it is again common to accelerate the preconditioner with an outer Krylov scheme.

From our brief discussion of these two model problems, it is clear that the design of linear solvers for multiphysics problems entails non-trivial choices:
If a block preconditioner is applied, as is usually the case, decisions are needed for a splitting strategy, subproblem solvers, and an outer Krylov method.
Each of these involves both categorical choices, i.e., which class of method to apply, and tuning of parameters within each solver.
Moreover, as the transport problem indicates, the optimal solver strategy may change dynamically with the solution.
For instance, for flow and transport problems, changes in well configuration can lead to frequent changes in the flow field, that may impact solver performance.
The strategy should also be adapted to changes in simulation models, for example varying geometry, boundary conditions etc.
As linear solvers often consume a substantial part of the overall computational time, these choices are of paramount importance.
While the importance of expert knowledge in selecting solvers relative to the mathematical structure of the problem is clear, users of simulation tools often lack direct access to such expertise. Consequently, they have to rely on heuristics for designing linear solvers.

The aim of this work is to develop an automated framework for solver selection. It seeks to identify optimal solvers, in terms of computational time, by integrating expert knowledge where available, with an exploration of solver options. Moreover, the framework should be able to adapt the solver to changes in the characteristics of the solution.
The framework should be designed  to be seamlessly embedded into the simulation process, incurring minimal costs in computational time and requiring minimal alterations to the simulation workflow and software.

We emphasize that we consider the problem of solver selection different, and downstream, from the problem of development of solvers:
Throughout this work, we assume that for any linear system considered, there is at least one candidate linear solver strategy, and our task is to make optimal choices within this strategy.
As new solvers are introduced our framework, this opens up the possibility of computational savings, but it also broadens the range of options that must be considered in the selection process.
Furthermore, the theoretical study of solvers will often lead to better understanding of available algorithms, such as permissible ranges for solver parameter values, that can be used to narrow the choices that must be explored by our framework.

The solver selection problem is a specific application of a more general algorithm selection problem \cite{RICE197665,Smith-Miles2009,BISCHL201641}, which also includes a broad range of other applications such as machine-learning algorithm parameter tuning, compiler optimization, Boolean satisfiability problems, and search algorithms, see e.g. \cite{Smith-Miles2009,BISCHL201641} for more information.
The general task is to find an appropriate algorithm for a given problem, based on characteristics of the problem and past performance data for the candidate algorithms. The algorithm selection methods have found various applications related to numerical simulations, including, but not limited to, those discussed below.

An established approach to solver selection is formulated as a classification problem, see, e.g., \cite{Jessup2016Lighthouse,PYTHIA,Dongarra2006SANS,Bhowmick2010ADTreeLinearSolvers,Eller2012}. Each problem in a training data set is solved exhaustively by all available solvers, and a machine learning algorithm learns to map a problem to its optimal solver. This mapping is then exploited for problems not observed in the training data set.
The classification approach requires that the performance data is available for all candidate solvers, however, this is not a feasible for problems with frequently changing dynamics, as illustrated by subsurface transport with changing well configurations.
Incomplete performance data can be dealt with by viewing the selection process as optimization under uncertainty.
A natural approach is to apply Bayesian optimization, wherein a probabilistic model is constructed for the target function and used to guide the solver selection.
Works in this direction include \cite{liu2021gptune,Roy2021Bliss}, which all apply a probabilistic representation based on Gaussian Processes (GP)  \cite{rasmussen2006gaussian}, or generalizations thereof.

Enhancements in solver selection can be attained by conceptualizing a simulation not as a singular process to be optimized, but as a series of successive problems. Each time step, or each iterative linear solution within a time step or within a nonlinear simulation, represents a distinct solver selection problem.
This entails a much deeper embedding of the solver selection process into the numerical simulation, which shortens the feedback loop and provides frequent information on solver performance.
The iterative nature of time stepping is exploited in \cite{Mishev2009} to select iterative solvers, preconditioners, and tune solver parameters. Additionally, the framework can exploit performance data in databases from previous simulations.
Selecting the linear solver for individual time steps also opens up the possibility of tuning the solver to changing characteristics of the solution,
as was done, e.g., in \cite{Bhowmick2006PETSc-FUN3D,Clees2010_alphaSAMG}. 
Common to these works is the manual definition of a problem-specific heuristic for when to switch solvers.
Finally, we note that these characteristics can be learned \cite{Silva2021}.

Our prime interest herein is to minimize the computational cost of solving linear system, hence our focus on performance data.
Other techniques can also be applied to obtain information on the system, and thus speed up simulations.
This includes the use of a posteriori error estimates within inexact Newton schemes, as was done in e.g. \cite{jiranek2010posteriori,ahmed2020adaptive}. 
We also note that there has recently been considerable interest in applying machine learning to tune parameters in multigrid methods, see e.g., \cite{caldana2023deep,RAMESHKUMAR2024112570,Huang2023}.

In this work, we build a solver selection framework which does not rely on a prior database of solver performance and also adapt to changes in the solution.
To that end, our framework relies on two key ingredients: 
First, we consider the selection of a solver for every time step as a separate optimization problem and update the selection model with newly obtained performance data. 
This facilitates automatic experimentation with solver configurations with a short feedback loop and thus mitigates the lack of comprehensive solver performance data prior to simulation.
Second, we build a relation between the performance of solver options and the characteristics of the solution and use this to steer the adaptation between time steps.
To the best of our knowledge, this is the first framework that combines these two features for the purpose of solver selection.
It should be noted that, while we limit solver selection to individual time-steps in the current work, the framework is also applicable to solver selection for iterative solutions of non-linear problems that are not time-dependent.

The rest of the manuscript is structured as follows:
Section \ref{sec:problem_statement} introduces the terminology needed to formally define the solver selection problem, while our methodology is presented in Section \ref{sec:solver_selection_framework}.
In Section \ref{sec:model_problems}, we present the governing equations for the two model problems of poroelasticity and non-isothermal coupled flow and transport, together with popular linear solver approaches for these problems.
The numerical examples presented in Section \ref{sec:numerical_experiments} show that our framework is indeed capable of selecting an optimal solver and to adapt the choice to changing solution characteristics. We also investigate different approaches to exploring uncertainty caused by scarce data, and show how the data from one simulation setup can be utilized for another simulation setup to speed up the selection process for it.
Finally, our concluding remarks are given in Section \ref{sec:conclusions}.

\section{Problem statement}
\label{sec:problem_statement}
In this section, we introduce the solver selection problem which consists of three components: The space of available solver configurations, a characterization of the static and dynamic features of the simulation, and the reward function.

\subsection{The solver configuration space}
\label{sec:solver_space}
\begin{figure}[htbp]
    \centering
    \begin{subfigure}[b]{0.49\textwidth}
        \centering
        \includegraphics[height=7cm]{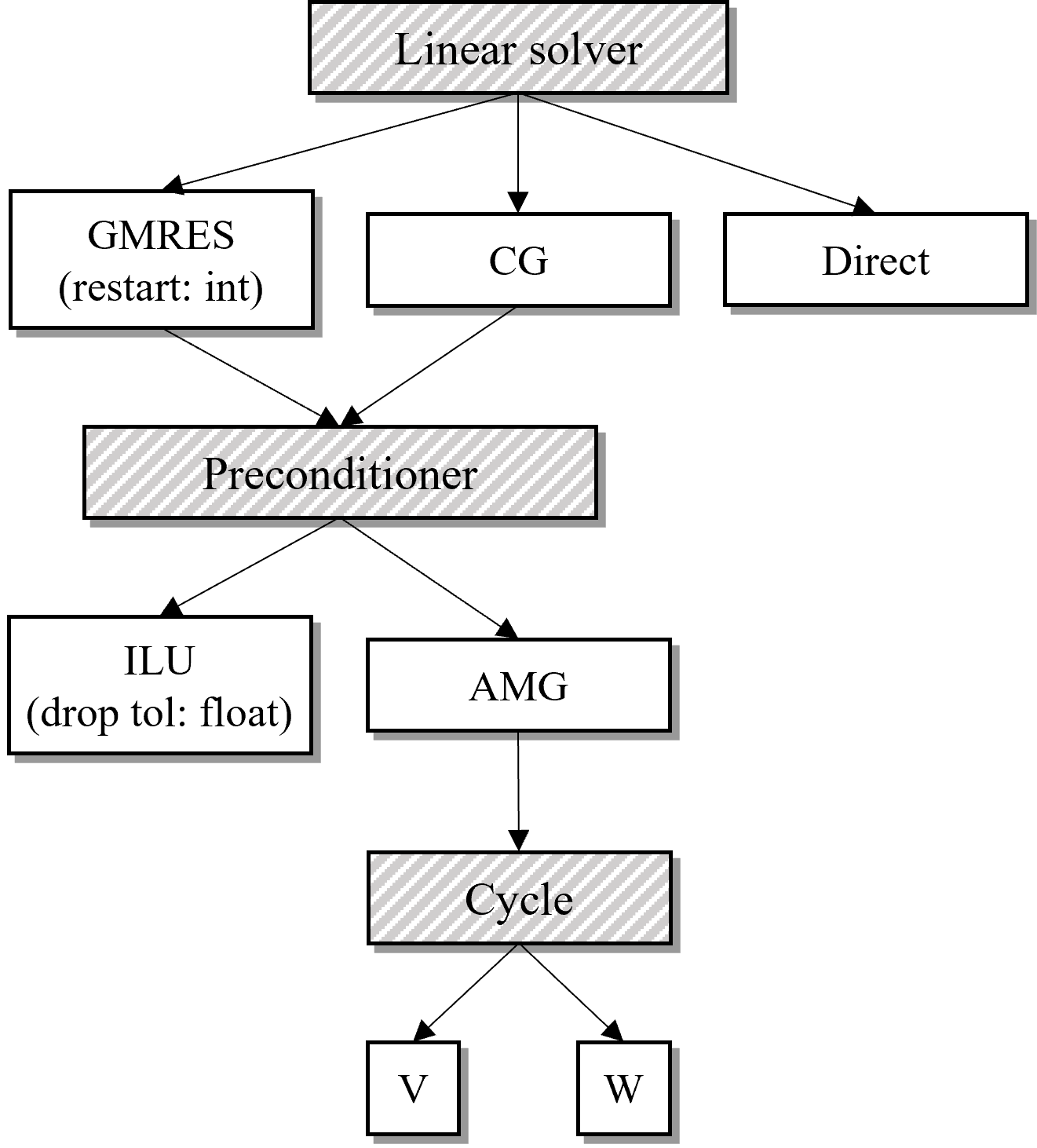}
        \caption{}
        \label{fig:solver_space}
    \end{subfigure}
    \begin{subfigure}[b]{0.49\textwidth}
        \centering
        \includegraphics[height=7cm]{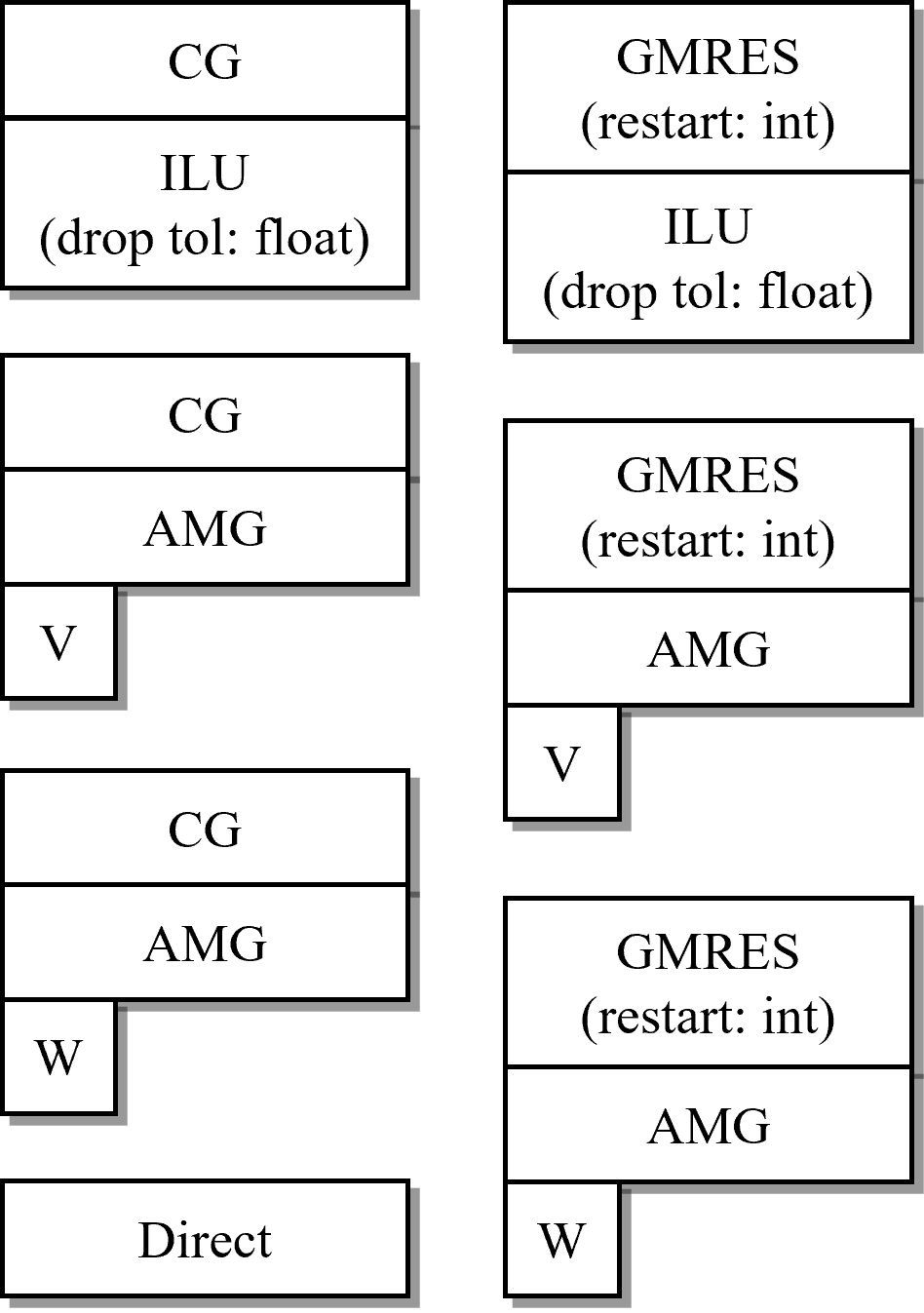}
        \caption{}
        \label{fig:independent_solvers}
    \end{subfigure}
    \caption{(a) Example solver space. The hatched blocks represent categorical decisions. (b) The seven solver configurations.}
\end{figure}

The solver selection framework seeks the optimal solver configuration in a user-provided solver configuration space which describes all candidate solver configurations.
The structure of this space can be non-trivial, as can be illustrated by an example for a linear single-physics problem, with a limited set of solver options:
The linear system is solved by either a direct method or an iterative Krylov solver, taken to be either the Generalized Minimal Residual (GMRES) or the Conjugate Gradient (CG) method.
If a Krylov solver is applied, it is preconditioned  by an Algebraic Multigrid (AMG) method or the incomplete LU (ILU) factorization.
Picking from these options leads to additional choices to be made, such as the restart parameter in GMRES, the drop tolerance in ILU and the type of cycle (V or W) in AMG, and so on. The solver configuration space for this example is visualized in Figure \ref{fig:solver_space}.

From this example, it is clear that deciding for a specific configuration involves two principal types of decisions: categorical and numerical. 
Categorical decisions entail choosing specific solver components, e.g., the Krylov scheme and the type of preconditioner, while numerical decisions involve the setting of real or integer values, e.g., the drop tolerance in ILU.
An important distinction between categorical and numerical decisions is that the latter has a natural ordering of its values which can be exploited by a solver selection algorithm.
Moreover, there is a finite number of combinations of categorical choices, while the numerical decisions that involve real numbers can take infinitely many values.
As can be seen in Figure \ref{fig:solver_space}, the choices have a structure reminiscent of a decision tree, in that a specific decision may trigger or exclude a later decision: The choice of preconditioner is only relevant if a Krylov solver was chosen, and the ILU drop tolerance is not needed in AMG.

We refer by the term \textbf{solver configuration} to the complete algorithm capable of solving a linear system. 
The set of all possible combinations of categorical choices is denoted $\mathcal{A}_\text{cat}$, and, for any $a_\text{cat}\in\mathcal{A}_\text{cat}$, we let $\mathcal{A}_\text{num}(a_\text{cat})$ denote the corresponding space of numerical choices.
The space of solver configurations,  $\mathcal{A}$, can be represented as
\begin{equation*}
    \mathcal{A} = \bigcup_{a_\text{cat}\in\mathcal{A}_\text{cat}} \mathcal{A}_\text{num}(a_\text{cat}).
\end{equation*}
In our example there are in total seven solver configuration groups, see Figure \ref{fig:independent_solvers}, but only some of the groups contain numerical parameters.

\subsection{The context space}
\label{sec:context}
The ideal solver configuration will depend on the properties of the problem at hand, such as the geometry of the model, the boundary and initial conditions, the source terms, the model parameters, etc.
Knowledge of such factors can be used to guide solver selection, as has also been shown previously \cite{Silva2021}.
In fact, this can be viewed as a mechanism for a user to provide the selection process with expert insight into the characteristics of a given simulation model.
To that end, a characterization which includes both static and dynamic features of the problem and solution, can be collected into an array of numerical values that we will refer to as the \textbf{simulation context}.

We define the simulation context as a space that includes all the possible values that the simulation context can take, $\mathcal{C} \subset \mathbb{R}^m$, where $m$ is the dimension of the context. 
We will provide the content for the specific problems in the corresponding experiments in Section \ref{sec:numerical_experiments}.

\subsection{Optimization problem}
\label{sec:optimizaton_problem}
The goal of the solver selection process is to find the optimal solver configuration $a^\star \in \mathcal{A}$
which minimizes the time spent on solving the linear system and penalizes failed attempts. 
The optimization problem is expressed in terms of the reward function $R: \mathcal{A} \times \mathcal{C} \rightarrow \mathbb{R}$.
Given the solver configuration $a\in \mathcal{A}$ and the context of the simulated model $c \in \mathcal{C}$, $R$ evaluates the solution efficiency:
\begin{equation}
\label{eq:reward}
    R(a; c) := 
    \begin{cases}
    -\log T_{\text{sol}}, & \text{if successfully solved}, \\
    P, & \text{otherwise}
\end{cases}
\end{equation}
where $T_\text{sol}$ is the time spent on solving the current problem. 
The penalty for a failed attempt, $P$, is set to twice the largest solver runtime for a given problem, that is:
\begin{equation*}
    P := -\log \left( 2\cdot \max{}T_{\text{sol}} \right)
\end{equation*}

The value of the reward function for a given solver configuration $a$ and context $c$ is not deterministic; therefore, $R$ should be considered a stochastic quantity.
There are in general two sources to the lack of determinism:
First, the context space may not fully describe the features of the simulation which should inform the solver selection;
indeed, if the context space is void, no such information is available.
Thus, for a specific configuration, $R$ may take multiple values for the same context. 
Although this can be alleviated by employing a rich description of the context, such an approach risks putting a heavy burden on the solver-selection method.
Second, the measurement of the runtime is associated with stochastic measurement noise caused by other programs executed on the same computer at the time of measurement.
This can in principle be corrected by carefully monitoring the processing running on the computer, but adjusting for such factors is technically challenging.
A different workaround, to construct the reward function based on iteration count instead of runtime, is not practical since candidate solver configurations will have different runtime cost for their iterations.
We chose instead to treat $R$ as a stochastic function and assume it is described by an unknown stationary distribution. 
We define the optimal solver configuration $a^\star$ by incorporating the stochastic nature of the reward function:
\begin{equation}
\label{eq:optimization_reward}
    a^\star = \argmax_{a \in \mathcal{A}} \mathbb{E} \left[ 
 R(a; c) \right]
\end{equation}
where $\mathbb{E} \left[ \cdot \right]$ denotes the expectation.
That is, we define the optimal solver as the one which produces the highest reward on average.
In the following, we will denote by \textbf{sampling the reward function} the process of applying a solver configuration to a problem and evaluating the corresponding reward.

\section{Solver selection framework}
\label{sec:solver_selection_framework}
Since the distribution of the objective function is unknown, we do not have a functional expression for its expectation. Thus, the formula Eq. \eqref{eq:optimization_reward} defines the optimal solver, but it does not provide a recipe for how to find it. 
In practice, we can only approximate the expectation of $R$ by applying black-box optimization methods and sampling $R$ with various solver configurations $a$ and contexts $c$.
As our interest is problems where each sampling is considered a computationally costly operation, it is critical that the optimization method efficiently deals with data scarcity. 

Our optimization approach consists of two main ingredients:
First, we apply online learning by assigning, in principle, a new solver for each time step in a simulation.
This provides a short feedback loop for the optimization tool and allows for rapid adjustment based on newly generated and problem-specific data from $R$.
For simplicity, the solver is fixed for all Newton iterations within one simulation time step. 
Second, we use machine learning regression methods to build a model for predicting $\mathbb{E}\left[ R(a; c) \right]$ and selecting solvers.
The solver selection will need to balance \textbf{exploitation}, where the best solver is selected based on currently available data, with \textbf{exploration}, which seeks to address lack of data in regions of the solver or context space.

\subsection{Online learning}
We consider each successive time-step as an independent problem to which we apply the solver selection algorithm.
All linear solutions within a time step generate performance data that is used to update the solver selection model during simulation, in what we will refer to as \textbf{online learning}.
This allows us to start a simulation without any performance data, but data from previous simulations can be included.
As the solver selection aims to minimize the time to solve each linear system, the process will also further the goal of minimizing the total execution time of the whole simulation.

Compared to a hypothetical setting where the ideal solver is known in advance, there are two sources of additional computational cost under online learning:
First, with limited performance data available, exploration of the solver configuration space will necessarily lead to suboptimal solvers being chosen. 
This cost can to some degree be alleviated by the selection algorithm, but it is also an unavoidable consequence of lack of information.
Second, time spent on updating the solver selection model with new performance data and selecting a solver should be considered part of the overall simulation time.
It is important to note that, while linear system solution time scales with the system size, the cost of the optimization algorithm will typically scale with the number of data points available.
Thus the cost of optimization will matter less for large problems where identifying appropriate solvers can be expected to bring the largest absolute decreases in simulation time.

\subsection{Machine learning methods}
\label{sec:machine_learning_methods}
To approximate $\mathbb{E}\left[ R(a; c) \right]$ we employ machine learning regression methods.
As discussed in Section \ref{sec:solver_space}, a solver configuration $a$ consists of categorical decisions of subalgorithms and their corresponding parameters.
Accordingly, we split the solver configuration into $a:=\left[a_\text{cat}, a_\text{num}\right]$.
For brevity, we denote by $x$ a pair of the solver configuration and the simulation context: $x := \left[ a, c \right] \in \mathcal{A} \times \mathcal{C} $; the reward function Eq. \eqref{eq:reward} can take this argument: $R(x) := R(a; c)$. 
The machine learning algorithms will rely on performance data, with the set of data points $D$ defined as:
\[
    D = \left\{ x_i, r_i \right\}_{i=1..n},
\]
where $r_i$ is the value of the reward function sampled for $x_i$, and $n$ is the number of data points in $D$. 
No distinction is made between data from previous simulations and the current simulation.

For the purpose of improving the accuracy of a machine learning algorithm, the data set is transformed to have zero mean and unit variance.
For numerical parameters that take real numbers, we split the admissible interval into a finite set of evenly distributed values and let the optimization choose between them. In practice, we use $20$ points for each numerical parameter.

Different elements in $\mathcal{A}_\text{cat}$ can share components, for instance both GMRES and CG can use the ILU preconditioner, referring to Figure \ref{fig:independent_solvers}.
Apparently, this would allow for sharing performance data between solver configurations and hence a faster reduction in the uncertainty of $R$. 
However, while it is easy to measure the time spent in the individual components of a solver, this information is highly dependent on the interaction with other solver components and will have limited to no value if transferred to a different set of categorical choices. 
Therefore, we disregard the possibility of such transfer, and instead assign different machine learning instances to each $a_\text{cat}\in\mathcal{A}_\text{cat}$, motivated by \cite{nguyen2020bayesian}.

If no data is available for a solver configuration at the start of the simulation, the first $n_{init}$ data points for each $a_{cat}$ are always collected via random exploration to provide the machine learning model with sufficient data to make meaningful predictions. In practice, we have found $n_{init}=3$ to be a reasonable parameter value.

For the machine learning method trained on data $D$, we denote by $\mathbb{E}_D[R(x)]$ the approximation of the expectation $\mathbb{E}[R(x)]$ defined in Section \ref{eq:optimization_reward}. 
In the following, we will outline two specific nonlinear regression methods which differ in the way they treat uncertainty, as the first incorporates Bayesian statistics, while the second applies heuristics.

\subsection{Bayesian optimization approach}
\label{sec:gp}
The classical technique for expensive black-box optimization under data scarcity is Bayesian optimization \cite{bayesian_optimization}. We give a brief overview of it here; for more information, see \cite{rasmussen2006gaussian,bayesian_optimization}. The core idea is to build a probabilistic approximation of the expected value and variance of the target function, based on prior data.
In practice, the Gaussian Process (GP) statistical model is applied.
The GP algorithm considers a random variable $R(x)$ for each point $x\in \mathcal{A} \times \mathcal{C} $ 
and assumes that for a finite number of points $x$, the values $R(x)$ follow a multivariate normal distribution.
The variance at a point $x$ given data $D$ is denoted $\sigma^2_D[R(x)]$ and will increase as we move away from available data, thus the GP representation allows us to explicitly estimate our confidence in the current model.

The decision-making process of the GP algorithm is guided by an acquisition function. We utilize the Upper Confidence Bound (UCB) acquisition function:
\begin{equation}
\label{eq:acquisition}
\text{UCB}(x) = \mathbb{E}_D[R(x)] + \alpha \cdot \sigma^2_D[R(x)],   
\end{equation}
where $\alpha$ is a parameter; we note that other acquisition functions are also used, for example, in \cite{Roy2021Bliss}. The next solver to be evaluated $\tilde{a}^\star$ is selected by maximizing the acquisition function: 
\[
\tilde{a}^\star = \argmax_{a \in \mathcal{A}} \text{UCB}(x).
\]
The chosen configuration $\tilde{a}^\star$ attempts to approach the true optimum $a^\star$ of Eq. \eqref{eq:optimization_reward}  by balancing exploration and exploitation: 
For low uncertainty, the first term will dominate the acquisition function and the solver with the highest expected reward is chosen, thus exploiting existing data.
If the second term in Eq. \eqref{eq:acquisition} dominates, the selection will explore parts of the solver configuration space with high uncertainty.
The parameter $\alpha$ determines how much preference to give to exploration. 
The space of functions predicted by the GP algorithm can be controlled by the GP kernel, which in this work is taken as 
the radial basis function kernel:
\[
k(x_i, x_j) = \exp{\left( -\dfrac{d(x_i, x_j)^2}{2l^2} \right)},
\]
where $d(\cdot, \cdot)$ is Euclidean distance and  $l>0$ is the algorithm parameter that corresponds to the length scale of the modeled function. This parameter is learned during fitting the GP algorithm with the data. 

As outlined in Section \ref{sec:machine_learning_methods}, we use separate instances of the GP method for different combinations of categorical decisions. Each instance is trained with the corresponding subset of data from $D$.
Selection of a new solver configuration entails evaluating all the GPs and picking the configuration that corresponds to the best prediction of all the GP instances.

GP is a powerful and established algorithm, which is commonly used in the solver selection domain \cite{liu2021gptune,Roy2021Bliss}.
Its computational cost can be significant, though: The cost of  updating the GP with new data is cubic with respect to the number of data points in $D$, and the cost of selecting a new solver is cubic with respect to the number of solver configurations to be selected from.
Hence using GP for solver selection may introduce a significant overhead, in particular for simulations that produce many data points or involve relatively small linear systems.
The performance of GP and the associated overhead will be tested in Section \ref{sec:numerical_experiments}.

\subsection{Heuristic approach}
\label{sec:eps_greedy}
As a computationally cheap alternative to GP, we here consider a heuristic-based approach to the black-box optimization.
This approach does not apply a representation of variance, thus we cannot explicitly estimate the uncertainty and utilize the UCB acquisition function.
Nevertheless, we can still approximate the expectation of the reward function by applying any machine learning regression model.  

To balance between exploration and exploitation we apply an epsilon-greedy heuristic.
With probability $\varepsilon \in [0, 1]$, we select a random solver configuration from $\mathcal{A}$. Otherwise, we make a greedy choice by selecting the configuration $\tilde{a}^\star$ that maximizes the current approximation of the expectation:
\[
\tilde{a}^\star = \argmax_{a \in \mathcal{A}} \mathbb{E}_D \left[ 
 R(x) \right]
\]
Note that $\tilde{a}^\star \neq a^\star$ from Eq. \eqref{eq:optimization_reward} because the former represents our current belief based on the data $D$, and the latter represents the true optimum, which is unknown.
After each exploration, the probability of exploring next time decreases: $\varepsilon \leftarrow \varepsilon \cdot \gamma$\, where $\gamma$ is a parameter.

The expectation of the reward function is approximated by the gradient boosting regression algorithm \cite{Friedman2001}. This machine learning model approximates data with piecewise constant functions and supports non-linear and discontinuous target functions.
Once fitted, the gradient boosting algorithm cannot be updated with newly generated data, so the model must be refitted from scratch after each new time step.
The cost of this operation is however insignificant, as will be shown in Section \ref{sec:numerical_experiments}.

As discussed in Section \ref{sec:machine_learning_methods}, separate instances of the gradient boosting algorithm are used to model the reward for different categorical decisions $a_\text{cat}$. The parameters $\varepsilon$, $\gamma$ and $n_{init}$ are separate for each instance, and 
if one algorithm instance decides to explore, it is selected regardless of how high the predicted rewards are from the other algorithm instances. Thus, the overall probability to explore is $1-(1-\varepsilon_1 )\cdot (1-\varepsilon_2 )\cdot \ldots \cdot (1-\varepsilon_k)=1-\Pi_{i=1}^k (1-\varepsilon_i)$ for $k$ instances of the algorithm.

\section{Model problems}
\label{sec:model_problems}
In this Section, we define the two model problems that will be used to test the solver selection framework: Deformation of a poroelastic medium and non-isothermal porous media flow of a single-phase fluid.
For each problem, we present the governing equations, discretization, and linear solver algorithms.
The following information is common for both models:
The porous medium is defined in the domain $\Omega \in \mathbb{R}^d$, where $d$ is the spatial dimension. The domain boundary consists of two disjoint subsets: $\partial \Omega = \Gamma_N \cup \Gamma_D$ that correspond to boundary conditions Neumann ($\Gamma_N$) and Dirichlet ($\Gamma_D$) boundary conditions. The problem is time-dependent with $t \in [0, t_\text{end}]$.

\subsection{Poromechanical deformation}
\label{sec:poromechanics}

The widely used mathematical model for poroelasticity is the quasi-static linear Biot model. The problem is based on mass conservation and force balance equations (see, e.g., \cite{coussy2004poromechanics}) and can be stated: Find the displacement $\mathbf{u}$ and the pressure $p$ such that

\begin{equation}
\label{eq:poromechanics}
\begin{cases}
    \nabla \cdot \sigma = \mathbf{0} 
    \\
    \rho \dfrac{\partial \phi}{\partial t} + \rho \nabla \cdot \mathbf{q} = f     
\end{cases}\quad \text{in } \Omega \times [0, t_\text{end}]
\end{equation}
Here $\rho$ is the fluid density, $\phi$ is the rock porosity, $\mathbf{q}$ is the fluid velocity, $f$ is the mass source term and $\sigma$ is the total stress. We apply Dirichlet and Neumann boundary conditions on the corresponding parts of the boundary $\Gamma_N$ and $\Gamma_D$. 

\subsubsection*{Constitutive laws}
The fluid velocity is governed by Darcy's law,
\begin{equation}
\label{eq:darcy_law}
    \mathbf{q} = -\dfrac{\mathbf{K}}{\mu} \nabla p,
\end{equation}
where $\mathbf{K}$ is the permeability tensor, $\mu$ is the fluid dynamic viscosity and gravitational effects are neglected. The total stress is defined as
\begin{equation}
    \mathbf{\sigma} = \mathbf{\sigma'} - bp\mathbf{I},
\end{equation}
where $\sigma'$ is the effective stress, $b$ is Biot's coefficient and $\mathbf{I}$ is the second-order unit tensor. The effective stress is connected to the displacement via Hooke's law (assuming zero reference stress and displacement) by
\begin{equation}
    \mathbf{\sigma'} = \mathbb{C} : \mathbf{\epsilon},
\end{equation}
where $\mathbb{C}$ is a fourth-order tensor stiffness and $\epsilon$ is defined as follows:
\begin{equation}
    \mathbf{\epsilon} = \nabla^s \mathbf{u} := \dfrac{1}{2} \left( \nabla \mathbf{u} + \nabla \mathbf{u}^T \right).
\end{equation}
We consider the isotropic homogeneous rock, so Hooke's law can be simplified to
\begin{equation}
    \sigma' = 2 \mu_r \epsilon(\mathbf{u}) + \lambda \nabla \cdot \mathbf{u} \mathbf{I},
\end{equation}
where $\mu_r$ and $\lambda$ are the Lame parameters.
In our poromechanics model, the fluid density and viscosity are constant, while porosity is a function of pressure and displacement,
\begin{equation}
    \phi(p, \mathbf{u}) = \dfrac{1}{m}p + b \nabla \cdot \mathbf{u},
\end{equation}
where $\dfrac{1}{m}$ is the specific storage.

\subsubsection*{Discretization}
We use fully implicit Euler time discretization. The spatial discretization is based on cell-centered finite volume schemes. We apply the two-point flux approximation \cite{Aziz1979} defined for diffusive scalar problems and the multi-point stress approximation (MPSA) for vector problems \cite{Nordbotten2016}. The spatial discretization is locally conservative for mass and momentum.

\subsubsection*{Solver configurations}
\label{sec:mandel_solvers}
For the given model problem, each time step of the simulation requires solving a linear system of equations.
We solve it with GMRES preconditioned by the Fixed-Stress approximation of the Schur complement (see \cite{White2016,Mikelic_Wheeler_fixed_stress}). The Fixed-Stress scheme is briefly reviewed below.

We can rewrite the system \eqref{eq:poromechanics} in the residual form by moving the source term to the left-hand side and denoting $r_\text{mass}$ and $r_\text{force}$ as the residuals of the mass and force balance equations, respectively. The corresponding matrix has the following block structure:
\begin{equation}
J = 
    \begin{bmatrix}
    A & B_1 \\
    B_2 & C \\
    \end{bmatrix},
\end{equation}
where $A = {\partial r_\text{force}}/{\partial \mathbf{u}}$, $B_1 = {\partial r_\text{force}}/{\partial p}$ , $B_2 = {\partial r_\text{mass}}/{\partial \mathbf{u}}$ and $C = {\partial r_\text{mass}}/{\partial p}$. The matrix can be factorized by a block $LDU$ decomposition:
\begin{equation}
    J = \begin{bmatrix}
        I & \\
        B_2 A^{-1} & I \\
    \end{bmatrix}
    \begin{bmatrix}
        A & \\
          & S_A
    \end{bmatrix}
    \begin{bmatrix}
        I & A^{-1}B_1 \\
          & I
    \end{bmatrix},
\end{equation}
where the term $S_A$ is the Schur complement with respect to $A$:
\begin{equation}
    S_A = C - B_2 A^{-1} B_1.
\end{equation}
We approximate $J$ by the block upper-triangular matrix: 
\begin{equation}
    J^{-1} \approx U^{-1} D^{-1} =
    \begin{bmatrix}
        A^{-1} & -A^{-1}B_1  S_A^{-1} \\
           & S_A^{-1}
    \end{bmatrix}.
\end{equation}
The Schur complement should not be computed explicitly as it requires inverting the matrix $A$, thus an approximation must be sought. 
The most common approximation for the poromechanics system \eqref{eq:poromechanics} is to the so-called Fixed-Stress preconditioner, which for finite volume discretizations is defined as follows:
\begin{equation}
    -B_2 A^{-1} B_1 \approx L \dfrac{\rho V}{\Delta t}.
\end{equation}
Here $V$ is the cell volume and $\Delta t$ is the simulation time step. The stabilization parameter $L\in \mathbb{R}$ has a significant impact on the convergence rate of the iterations.
In a recent study \cite{Storvik2018}, the analytical limits for the optimal value of $L$ were derived: $L_\text{opt} \in \left[ L_\text{min}, L_\text{phys} \right]$, where
\begin{align}
    L_\text{min} = \dfrac{b^2}{4\mu_r + 2\lambda} \text{,}\quad
    L_\text{phys} = \dfrac{b^2}{\dfrac{2\mu_r}{d} + \lambda}.
\end{align}
However, it was shown in \cite{Both2019}, that for a given problem, the value of $L_\text{opt}$ depends on the boundary conditions and flow parameters.

\subsection{Non-isothermal flow in porous media}
\label{sec:thermal}
Assuming single-phase flow in porous media, a model that combines flow and thermal transport is based on mass and energy conservation equations, which read, respectively:
\begin{equation}
\label{eq:thermal_system}
\begin{cases}
    \phi \dfrac{\partial \rho }{\partial t} + \nabla \cdot (\rho \mathbf{q}) = f \\[10pt]
    \phi \dfrac{\partial h}{\partial t} + (1 - \phi) \dfrac{\partial h_r}{\partial t} + \nabla \cdot (h \mathbf{q}) + \nabla \cdot \mathbf{q_T} = f_T
\end{cases}\quad \text{in }  \Omega \times [0, t_\text{end}]
\end{equation}
where $h$ is the specific enthalpy of the fluid, $h_r$ is the specific enthalpy of the rock, $\mathbf{q_T}$ is the diffusive heat flux and $f_T$ is the energy source. The primary variables are pressure $p$ and temperature $T$. We apply Dirichlet and Neumann boundary conditions on the corresponding parts of the boundary $\Gamma_N$ and $\Gamma_D$.

We use the Peaceman well model \cite{Peaceman1978,Chen2009} to simulate wells; details can be found in \ref{appendix:thermal}.

\subsubsection*{Constitutive laws}

In the non-isothermal flow model in porous media, porosity is constant and gravity is neglected. The density is a function of pressure and temperature while the viscosity is a function of temperature; the details are given in \ref{appendix:thermal}. The specific enthalpy of fluid and rock are defined by the following simplified expressions: $h = \rho c_v T$, $h_r = \rho_r c_r T$, where $\rho_r$ is the density of rock and $c_v$ and $c_r$ are the specific heat capacities of fluid and rock, respectively.
The diffusive heat flux is defined by Fourier's law:
$\mathbf{q_T} = -k_T \cdot \nabla T$, where $k_T=\phi k_{Tf} + (1 - \phi) k_{Tr}$; $k_{Tf}$ is a thermal conductivity of fluid, and $k_{Tr}$ is a thermal conductivity of rock.

\subsubsection*{Discretization}
We use fully implicit Euler time discretization and finite volume methods for spatial derivative: Diffusive terms are discretized by two-point flux approximation, while advective fluxes are discretized by a first-order upstream scheme.

\subsubsection*{Solver configurations}
\label{sec:thermal_solvers}
The non-isothermal flow model is nonlinear and is solved by Newton's method. We adjust the time step dynamically so that Newton's method converges in around 4 iterations. 
At each Newton iteration, we solve the full Jacobian linear system with GMRES. Restricted to each time step, 
we use as preconditioner either the Constrained Pressure Residual (CPR) \cite{Wallis1983,Wallis1985} or the or the Schur complement approximation \cite{Roy2019}, both of which are briefly outlined below.

Similar to the poromechanics system from Section \ref{sec:poromechanics}, we rewrite Eq. \eqref{eq:thermal_system} in residual form, and define $r_\text{mass}$ and $r_\text{energy}$ as the residuals for the mass and energy equations, respectively.
Borrowing notation from \cite{Roy2019}, the block structure of the Jacobian matrix can be written 
\begin{equation}
\label{eq:thermal_jacobian}
J = 
    \begin{bmatrix}
        A_{pp} && A_{pT} \\
        A_{Tp} && A_{TT} \\
    \end{bmatrix},
\end{equation}
where $A_{pp} = {\partial r_\text{mass}}/{\partial p}$, $A_{pT} = {\partial r_\text{mass}}/{\partial T}$, $A_{Tp} = {\partial r_\text{energy}}/{\partial p}$ and $A_{TT} = {\partial r_\text{energy}}/{\partial T}$.

The CPR preconditioner is a multiplicative two-step preconditioner. The algorithm consists of the application of two operators: $M_1^{-1}$ and $M_2^{-1}$. For a right-hand side vector $b$, it is defined as follows:
\begin{enumerate}
    \item Precondition using $M_1: x_1 = M_1^{-1} \cdot b$
    \item Compute the new residual: $b_1 = b - J\cdot x_1$
    \item Precondition using $M_2$ and correct:
    $x = M_2^{-1}\cdot b_1 + x_1$
\end{enumerate}
The operator $M_1$ approximately inverts the pressure block:
\begin{equation}
    M_1^{-1} \approx \begin{bmatrix}
        A_{pp}^{-1} & 0 \\
        0 & 0
    \end{bmatrix},
\end{equation}
where $A_{pp}^{-1}$ is approximated using one AMG V-cycle. The operator $M_2$ approximately inverts the full system $J$ using an incomplete LU factorization (ILU).
While the CPR preconditioner applied to multi-phase flow usually includes a decoupling operator, e.g. \cite{Lacroix2003}, this was shown not to lead to additional improvement for the non-isothermal single-phase flow model \cite{Roy2019}, and we therefore do not apply any decoupling operator.

The second preconditioner is based on a Schur complement approximation. Again, we apply the block $LDU$ decomposition to the matrix \eqref{eq:thermal_jacobian}:
\begin{equation}
    J = \begin{bmatrix}
        I & 0 \\
        A_{Tp} A_{pp}^{-1} & I\\
    \end{bmatrix} \begin{bmatrix}
        A_{pp} & 0 \\
        0 & S_{T} \\
    \end{bmatrix} \begin{bmatrix}
        I & A_{pp}^{-1} A_{pT} \\
        0 & I
    \end{bmatrix},
\end{equation}
where $S_T = A_{TT} - A_{Tp} A_{pp}^{-1} A_{pT}$ is the Schur complement. The inverse of the block decomposition is given by:
\begin{equation}
\label{eq:thermal_jacobian_approximation}
    J^{-1} = \begin{bmatrix}
        I & -A_{pp}^{-1} A_{pT} \\
        0 & I
    \end{bmatrix} \begin{bmatrix}
        A_{pp}^{-1} & 0 \\
        0 & S_T^{-1} \\
    \end{bmatrix} \begin{bmatrix}
        I & 0 \\
        -A_{Tp} A_{pp}^{-1} & I
    \end{bmatrix}.
\end{equation}
Three different preconditioners can be derived by different approximations to this block inverse: What we will term the full preconditioners include all three factors, while upper and lower block triangular versions are found by including, respectively, the two first and the two last matrices.
We approximate the Schur complement as,
\begin{equation}
\label{eq:thermal_schur_approximation}
    \widetilde{S}_T = \phi \dfrac{c_v \rho}{\Delta t} + (1 - \phi) \dfrac{\rho_r c_r}{\Delta t} + \nabla \cdot (c_v \rho \mathbf{q}) - \nabla \cdot (k_T \nabla ) - c_v f_\text{prod},
\end{equation}
where $\Delta t$ is the simulation time step and $f_\text{prod}$ is the source term corresponding only to the production wells, see \cite{Roy2019} for a derivation. The preconditioner of Eq. \eqref{eq:thermal_jacobian_approximation} approximates each of the inverses of the matrices $A_{pp}$ and $\widetilde{S}_T$ by one AMG V-cycle iteration.

The quality of the approximation of $\widetilde{S}_T^{-1}$ depends on the ratio between the diffusive and advective terms in Eq. \eqref{eq:thermal_schur_approximation}, that is, of the Peclet number:
The AMG method can be expected to perform well for low Peclet numbers, thus diffusion dominates. In the advection-dominated regime, different approximations to $\widetilde{S}_T^{-1}$ are preferrable.
In contrast, CPR is designed for advection-dominated transport, and will likely suffer for low Peclet numbers.
The non-isothermal flow model problem is thus an example where the optimal solution strategy depends on the solution regime.
\section{Numerical experiments}
\label{sec:numerical_experiments}
To investigate the performance of the solver selection algorithm, we consider two test cases:
The Mandel problem in poromechanics, based on the model presented in Section \ref{sec:poromechanics}, and a case of non-isothermal flow and transport, as described in Section \ref{sec:thermal}. Full descriptions of the modeled systems are given in \ref{appendix:poro} and \ref{appendix:thermal}, respectively.

For these model problems we consider a series of experiments. In Section \ref{+sec:numerics:categorical}, we consider cases with 
few candidate solver configurations to choose from, but their optimal parameters must still be determined.
In the poromechanical simulation, we fix the categorical choices and focus only on optimizing the numerical value of the stabilization parameter $L$ described in Section \ref{sec:mandel_solvers}.
For the flow and transport problem, we seek when to switch between the CPR and the Schur-complement based preconditioners from Section \ref{sec:thermal_solvers}. This problem is set up with a cyclic injection pattern, which repeatedly brings the problem back to similar states, thus probing the selection algorithm's ability to exploit previous data. The cycle pattern is also constructed, by a somewhat contrived design of the injection rate and schedule, to make the system exhibit a wide range of Peclet numbers through the simulation, spanning regimes where both the CPR and Schur-complement based preconditioner have their forte.
In Section \ref{sec:experiment:extended_solver_space}, we consider the non-isothermal flow problem based on the same setup as before but extend the number of candidate solver configurations. This will test the selection algorithm's ability to efficiently zoom in on the best solvers.

As previously indicated, there are already well-established solver strategies for both model problems, and we do not expect to arrive at results that are significantly different from this existing knowledge. Indeed, our motivation for considering these cases is, in part, to study how efficiently the solver selection algorithm settles on strategies that could have been surmised from existing understanding of these problems but without relying on heuristics and manual tuning.

The solver selection framework is implemented in Python and is fully open-source\footnote{\href{https://github.com/Yuriyzabegaev/solver_selector}{github.com/Yuriyzabegaev/solver\_selector}}.
Moreover, we provide the Docker image with the reproducible experiments from this article \cite{docker_image}.
The simulation models are implemented in the open-source package PorePy \cite{Keilegavlen2021}, while the gradient boosting and Gaussian Process algorithms are provided by Scikit-learn \cite{scikit-learn}.
We use the implementation of GMRES and a direct solver provided by SciPy \cite{2020SciPy-NMeth}, as well as the flexible GMRES (FGMRES) implementation in PyAMG \cite{pyamg2023}.
As preconditioners, we use the ILU(0) algorithm provided by PETSc \cite{petsc-user-ref,DalcinPazKlerCosimo2011} and the Boomer AMG library \cite{boomeramg}. While other implementations may be more efficient, our motivation is to select among the provided options; thus, the important factor is to give the same settings to all options.

Linear systems are solved by GMRES with a restart after 30 inner iterations and with a maximum of 10 restarts. The linear solver relative tolerance is always set to $10^{-10}$.
For the gradient boosting algorithm and Boomer AMG, we use default parameters; parameters used for the Bayesian and heuristic optimization methods are discussed in Section \ref{sec:solver_selection_framework}.

To limit measurement noise, we run our experiments with no other CPU- or memory-intensive processes being executed on the same machine and use a single core.
We present the average and variance of 20 repetitions of all experiments, with the solver selection potentially taking different random choices in individual experiments.
Data is not shared between these repeated experiments.

\subsection{Proving the concept of automated solver selection}
\label{+sec:numerics:categorical}
In this set of experiments, we provide a small space of solver configurations to the solver selector, and show that solver selection improves simulation runtime compared to human-selected benchmarks. 
The optimization problem is solved using the heuristic approach defined in Section \ref{sec:eps_greedy} with the parameters $\varepsilon=0.5$ and $\gamma=0.9$. %

\subsubsection{Optimization of numerical solver configuration parameter}
To study selection of a numerical solver configuration parameter, we consider Mandel's problem and optimize the $L$ stabilization parameter of the Fixed-Stress preconditioner. 
The solver space for this problem is:
\[
\mathcal{A} = \left\{\text{GMRES} - \text{Fixed-Stress} - L \in \left[ L_\text{min}, L_\text{phys} \right] - \text{Direct subsolvers} \right\}.
\]
The optimal value of $L$ for this specific simulated model is believed to be static, so the context space $\mathcal{C}$ is void.

\begin{figure}[tbh]
    \centering
    \begin{subfigure}[b]{0.49\textwidth}
        \centering
        \includegraphics[width=\textwidth]{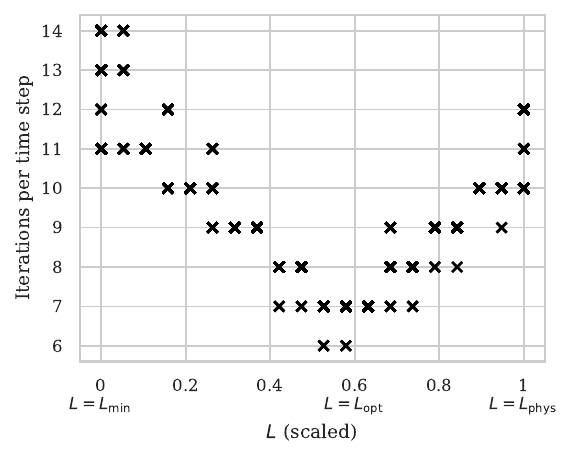}
        \caption{}
        \label{fig:poro_iters}
    \end{subfigure}
    \hfill
    \begin{subfigure}[b]{0.49\textwidth}
        \centering
        \includegraphics[width=\textwidth]{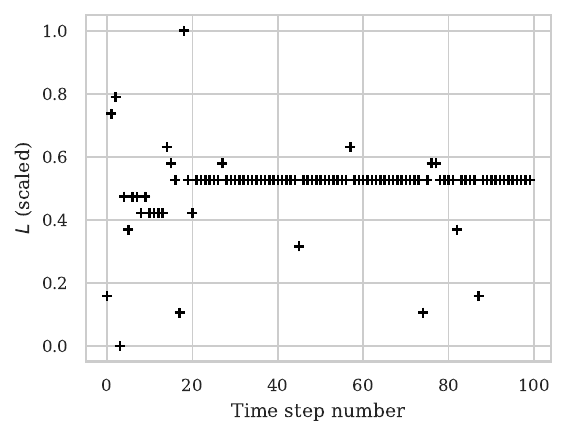}
        \caption{}
        \label{fig:poro_lfactor_dt}
    \end{subfigure}
    \caption{Optimization of numerical solver configuration parameter for the Mandel's problem. (a) Variation of number of GMRES iterations with the stabilization coefficient $L$. (b) Evolution in $L$ from the solver selection framework during one simulation.}
\end{figure}

Figure \ref{fig:poro_iters} shows how the number of GMRES iterations depends on the value of $L$ employed for a time step.
We apply a linear scaling to $L$ so that $L_\text{min}=0$ and $L_\text{phys} = 1$, and the values from all the time steps are shown together, so the results for the same $L$ can slightly vary from one time step to another.
From the figure, it is clear that the optimal value for this problem is approximately 0.6.
The search for the optimal value of $L$ for one of the 20 runs of this experiment is shown in Figure \ref{fig:poro_lfactor_dt}:
The exploration process starts with random sampling and evolves to primarily selecting the optimal value while reducing random decisions, following the epsilon-greedy exploration policy.
Continued exploration after the optimal value has been found, and the associated overhead in computational cost, can be seen as the price for safeguarding against dynamic changes in the optimal value.

\begin{figure}[tbh]
\centering
\begin{subfigure}{.5\textwidth}
  \centering
        \includegraphics[width=\linewidth]{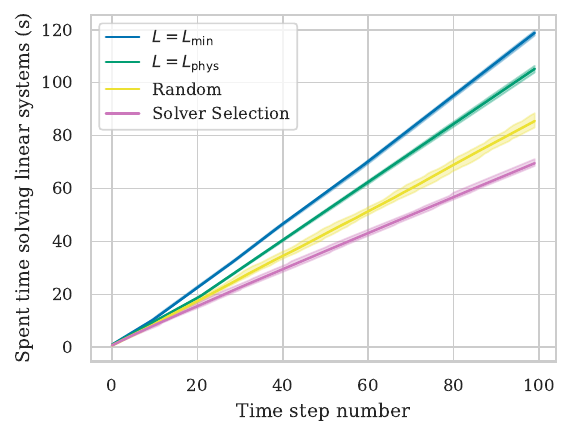}
        \caption{Mandel's problem}
        \label{fig:poro_performance}
\end{subfigure}%
\begin{subfigure}{.5\textwidth}
  \centering
        \includegraphics[width=\linewidth]{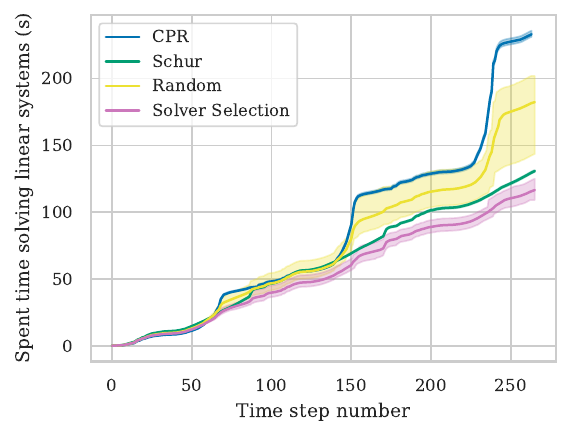}
        \caption{Non-isothermal flow problem}
        \label{fig:thermal_performance}
\end{subfigure}
\caption{Cumulative performance of solver configurations picked by the solver selector compared to pre-selected configurations and random decisions. Lower is better. The bold line of each color displays the mean value, while the lighter-shaded region of the same color around it indicates the range from the lowest to highest value among all the repeats of an experiment. }
\end{figure}

To measure the overall computational gains from applying the solver selection process, we consider the cumulative time spent on solving linear systems.
In Figure \ref{fig:poro_performance}, this is shown for four solver options: The solver selection algorithm, $L$ fixed to, respectively $L_\text{min}$ and $L_\text{max}$, and $L$ randomly chosen in each time step in the interval $[L_\text{min},L_\text{max}]$. As can be seen, the solver selection delivers the best results, while the relative success of the randomly chosen parameter can be contributed to the convex shape of the performance data shown in Figure \ref{fig:poro_iters}.
The latter point indicates that the test problem is relatively easy, but we reemphasize that the selection algorithm found the optimal $L$ value with no need for manual tuning, and with little overhead compared to running the full simulation with the optimal value of $L$.

\subsubsection{Optimization of categorical solver configuration choices}
\label{sec:numerics:categorical}
To study the impact of categorical decisions, we consider the non-isothermal flow problem, with the following two solver configurations:
\[
\mathcal{A} = \left\{ \left( \text{GMRES} - \text{CPR} \right); \left( \text{GMRES} - \text{Schur} \right) \right\}.
\]
All the numerical parameters are fixed. The Schur complement preconditioner in this experiment applies the full factorization of the block inverse Eq. \eqref{eq:thermal_jacobian_approximation}.
As noted in Section \ref{sec:thermal}, the optimal solver configuration is expected to depend on the dominant transport mechanism, measured by the Peclet number.
Accordingly, the Peclet number is included in the 
context space $\mathcal{C}$ together with other information listed in Table  \ref{tab:context_thermal}.
\begin{table}[tbh]
    \centering
    \renewcommand{\arraystretch}{1.1} %

    \begin{tabular}{l l}
        \toprule
        Parameter & Transformations \\
        \midrule
        Time step value, $\text{s}$ & $\log(\cdot)$ \\
        Peclet number & $\log \text{max} (\cdot)$ \\
        Peclet number & $\log \text{mean} (\cdot)$ \\
        Injection source rate, $\text{m}^3\text{s}^{-1}$ & $\log (\cdot)$\\
        Production source rate, $\text{m}^3\text{s}^{-1}$ & $\log (\cdot)$\\
        Well activity binary indicator& \\
        \bottomrule
    \end{tabular}
    \caption{The non-isothermal flow problem context space. The right column describes the transformations applied to these parameters.}
    \label{tab:context_thermal}
\end{table}

\begin{figure}[tbh]
    \centering
    \begin{subfigure}[b]{0.49\textwidth}
        \centering
        \includegraphics[width=\textwidth]{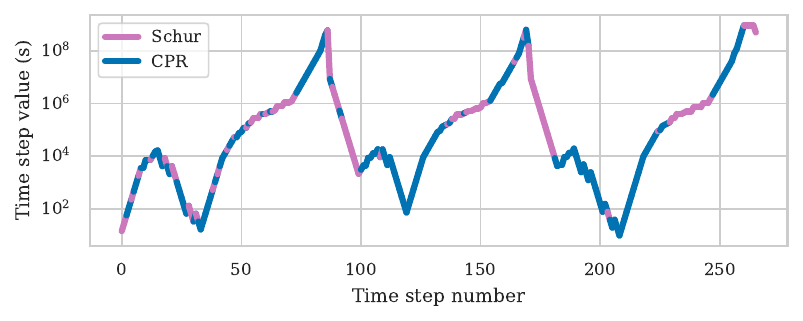}
        \caption{Evolution of time step magnitude.}
        \label{fig:thermal_dt}
    \end{subfigure}
    \hfill
    \begin{subfigure}[b]{0.49\textwidth}
        \centering
        \includegraphics[width=\textwidth]{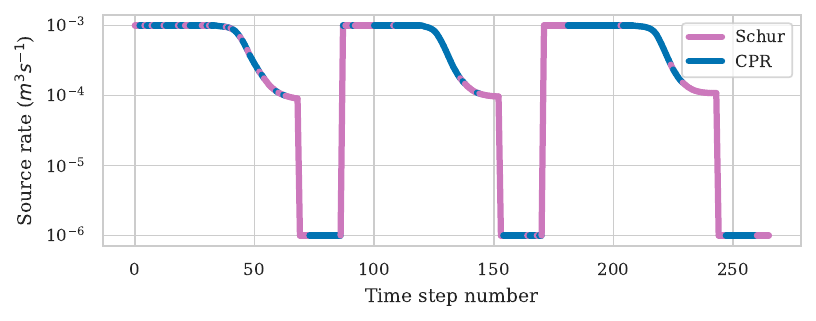}
        \caption{Evolution of inlet source rate magnitude.}
        \label{fig:thermal_source}
    \end{subfigure}
\caption{Solver selection for the non-isothermal flow problem. Different colors represent the solver which was used for each time step.}
\label{fig:thermal_solver_switching}
\end{figure}

The simulation goes through three injection periods that change the simulation characteristics. 
We can track the process of solver selection in Figure \ref{fig:thermal_solver_switching}, which shows the dependence of the time step magnitude and the volumetric injection rate on the time step number.
The time step size increases rapidly in periods without injection, followed by rapid decrease when injection is switched on.
If the Newton solver reaches the maximum number of iterations, we decrease the time step and select a solver configuration again. The performance data of these unsuccessful iterations is also used to update the machine learning algorithm, and they are also included in our reporting on linear solver behavior.

Figures \ref{fig:thermal_solver_switching} further show the choice of solver configuration for each time step. During the first injection period, the solver selector prioritizes exploration, with frequent switching between the CPR and Schur complement options. After approximately time step 30, the solver has collected sufficient data to settle on the following policy: CPR is preferred when the injection rate is high or the injection is switched off, while the Schur method is chosen for low injection rates. Periods of high and low injection rates correspond to high and low Peclet numbers, thus the transition from CPR to Schur is as expected, and in agreement with the results from \cite{Roy2019}. When injection is switched off, the system is diffusion dominated, but it is also close to steady state, thus the linear system is easier to solve. In this regime, the ILU method applied in CPR is sufficiently efficient to result in a lower total computational cost than that of the Schur method, although upon inspection, the difference turns out to be minor.

To evaluate the effect of solver selection for this problem, we compare the performance of the simulation driven by solver selection with three benchmarks:
Two preconditioners that use, respectively, the CPR and Schur complement preconditioners throughout the simulation, and one which randomly picks among them for each time step.
The cumulative time to solve linear systems is shown in Figure \ref{fig:thermal_performance}. 
For this specific problem setup, the relative cost of misusing CPR is high, while the relative cost of misusing Schur is relatively low. 
The extra cost of CPR is to a large degree associated with three short time periods which coincide with the injection well switching from rate to pressure control. 
The solver selector switches to the Schur complement method at these points, thus avoiding the problem.
When also accounting for the overhead associated with exploration, the solver selector brings only marginal gains compared to the Schur complement preconditioner.
However, solver performance will depend on the problem at hand, thus the optimal solver cannot be known \textit{a priori}. Moreover, the example also shows how the optimal solver may change during a simulation, and this effect is expected to be stronger if the problem complexity or range of solver options is increased.
The solver selector's ability to adapt to such changes gives it the potential to significantly reduce the computational burden of such simulations.

\subsubsection{Exploration sensitivity}
To better understand the behavior of the heuristic optimization algorithm, we investigate its sensitivity to the exploration parameter $\varepsilon$.
For both the poroelastic and non-isothermal flow setups, we repeated the experiments with a default value $\varepsilon = 0.5$ (same as in the previous experiment), with a higher value $\varepsilon=0.7$ and with zero exploration, $\varepsilon=0$. The latter case corresponds to a simple greedy algorithm which does not actively seek new data points, but may still obtain new data as part of the optimization loop.

The results are shown in Table \ref{tab:time_spent_exploration}, where we also provide the results of the random choice policy as a reference.
For the poroelastic model, neglecting exploration can significantly increase the spread of the results, as the method may or may not approach the minimum of the convex function to be optimized, cf. Figure \ref{fig:poro_iters}. Increased exploration gives insignificant changes to the standard model.
The structure of the non-isothermal flow problem is different, as its dynamic nature means the simulation characteristics change for every time step. Hence, during the first injection period, 
the machine learning algorithm always finds itself in an unexplored region of the space $\mathcal{C}$.
Predictions by machine learning algorithm outside the training data region are unreliable and chaotic, so it naturally leads to a forced exploration. Additional algorithmic exploration triggered by the $\varepsilon$ parameter prevents the worst-case result, and its overhead is negligible compared to the mean- and best-case results of the greedy algorithm. However, the amplified exploration leads to a worse result compared to the standard algorithm because, for this simulated model, the cost of suboptimal decisions is higher.

While these examples indicate that tuning of optimization parameters can give some computational savings, they also show that the impact of such changes is highly problem dependent. 
The main conclusion to be drawn from Table \ref{tab:time_spent_exploration} is that, while too much exploration in some cases can be considered an unnecessary overhead, there may also be considerable costs in not exploring.

\begin{table}[tbh]
    \centering    

\begin{tabular}{l | rrr | rrr}
\toprule
 & \multicolumn{3}{r}{Mandel's problem} & \multicolumn{3}{r}{Non-isothermal flow problem} \\
 & Worst & Mean & Best & Worst & Mean & Best \\
\midrule
Standard exploration $\varepsilon = 0.5$ & 71.68 & \textbf{69.40} & 67.61 & \textbf{122.75} & 114.87 & 107.93 \\
No exploration $\varepsilon = 0$ & 102.42 & 70.19 & \textbf{65.65} & 134.64 & \textbf{114.57} & \textbf{107.83} \\
More exploration $\varepsilon = 0.7$ & \textbf{71.30} & 69.73 & 68.17 & 135.63 & 117.41 & 108.23 \\
Random choice & 88.47 & 85.48 & 83.03 & 202.02 & 182.64 & 143.52 \\
\bottomrule
\end{tabular}

    \caption{Cumulative time (seconds) to solve linear systems for the whole simulation for the Mandel's problem and the non-isothermal flow problem. Results for the solver selection policy with different values of the $\varepsilon$ parameter and the result for the random policy. The minimal value in each column is given in bold.}
    \label{tab:time_spent_exploration}
\end{table}

\subsubsection{Knowledge transfer}
\label{sec:experiment:knowledge_transfer}
The simulation context allows for the transfer of performance data between simulations with different setups.
For this to be meaningful, the relative performance of different solvers must be similar in the setups, and this similarity must be seen by the simulation context. Borrowing from machine learning terminology, we refer to such knowledge transfer as warm-start simulations, in contrast to the cold-start simulations where all performance data must be collected through online learning.

To illustrate potential gains of warm starts, we alter the non-isothermal flow setup by moving the two wells to regions of lower permeability, see Figure \ref{fig:well_locations_permeability} in \ref{appendix:thermal}.
For both warm and cold start we consider 20 simulations to chart the impact of random variations. 
The warm-start simulations are given access to the data set $D$ from one of the experiments reported above.
The results are reported in Table \ref{tab:time_spent_knowledge_transfer}, where we note that, since we are solving a problem with different well locations, the time measurements cannot be directly compared to those listed in Table \ref{tab:time_spent_exploration}.
The warm start consistently leads to slightly better results. 
While this may be considered a modest gain, we recall from Figure \ref{fig:thermal_dt} that the initial phase of exploration, which is where the access to previous data can make a difference compared to the cold start, only lasted about 30 time steps with the original well locations and is likely of similar length in this case.
Again, this effect will grow with the problem complexity or number of solver options.

\begin{table}[tbh]
    \centering    
\begin{tabular}{l | rrr}
\toprule
 & Worst & Mean & Best \\
\midrule
Cold start & 59.17 & 55.42 & 50.14 \\
Warm start & \textbf{56.99} & \textbf{54.19} & \textbf{50.06} \\
\bottomrule
\end{tabular}
    
    \caption{Cumulative time (seconds) to solve linear systems for the whole simulation for the non-isothermal flow problem with the changed well locations. A minimal value in each column is in bold. 
   The time measurements cannot be directly compared to those listed in Table \ref{tab:time_spent_exploration}  because we are solving a problem with different well locations.}
    \label{tab:time_spent_knowledge_transfer}
\end{table}

\subsection{Extension of the available solver configuration space}
\label{sec:experiment:extended_solver_space}
Having studied the behavior of the solver selector applying the heuristic optimization on a very limited set of choices, we now introduce two extensions, focusing on the non-isothermal flow setup in Section \ref{sec:numerics:categorical}.
First, we consider a wider set of categorical choices with in total 19 different sets of options.
This expansion will in effect lengthen the exploration phase, in which the solver selector must deal with substantial uncertainty in the simulation data.
Our second extension is therefore to consider optimization based on Gaussian Processes in addition to the heuristic approach, with the motivation that GP's explicit representation of uncertainty may make a difference in the rate at which a policy is settled upon.

\begin{table}[tbh]
\centering
\begin{tabular}{r | lllll}
\toprule
Id & Outer solver & Preconditioner & Factorization & Thermal & Flow \\
\midrule
1 & GMRES & Schur & Lower & AMG & AMG \\
2 & GMRES & Schur & Upper & AMG & AMG \\
3 & GMRES & Schur & Full & AMG & AMG \\
4 & GMRES & CPR &  &  & AMG \\
5 & FGMRES & Schur & Lower & GMRES + AMG & GMRES + AMG \\
6 & FGMRES & Schur & Lower & AMG & GMRES + AMG \\
7 & FGMRES & Schur & Lower & GMRES + AMG & AMG \\
8 & FGMRES & Schur & Lower & AMG & AMG \\
9 & FGMRES & Schur & Upper & GMRES + AMG & GMRES + AMG \\
10 & FGMRES & Schur & Upper & AMG & GMRES + AMG \\
11 & FGMRES & Schur & Upper & GMRES + AMG & AMG \\
12 & FGMRES & Schur & Upper & AMG & AMG \\
13 & FGMRES & Schur & Full & GMRES + AMG & GMRES + AMG \\
14 & FGMRES & Schur & Full & AMG & GMRES + AMG \\
15 & FGMRES & Schur & Full & GMRES + AMG & AMG \\
16 & FGMRES & Schur & Full & AMG & AMG \\
17 & FGMRES & CPR &  &  & AMG \\
18 & FGMRES & CPR &  &  & GMRES + AMG \\
19 & Direct &  &  &  &  \\
\bottomrule
\end{tabular}
\caption{The solver configurations considered in Section \ref{sec:experiment:extended_solver_space}. Items that are not relevant for a given configuration are blank.}
\label{tab:extende_solver_options}
\end{table}

The different solver configurations are listed in Table \ref{tab:extende_solver_options} and can be summarized as follows:
In addition to solving the diffusive subproblems with AMG, as was done in Section \ref{sec:numerics:categorical}, we also consider AMG preconditioned with GMRES as inner solvers. This requires changing the outer Krylov solver to flexible GMRES (FGMRES).
For the preconditioners considered in Section \ref{sec:numerics:categorical}, i.e., without an inner GMRES method, we consider both GMRES and FGMRES as the outer Krylov method. Recall that, while these methods are identical in terms of the number of iterations, we consider implementations from two different libraries. Thus, this is an example of how the solver selector can deal with choices relating to software.
For the Schur complement method, we consider the three approximations of the inverse Jacobian Eq. \eqref{eq:thermal_jacobian_approximation} noted in Section \ref{sec:thermal}.
Finally, we add the option of solving the entire system with a direct solver.

\begin{figure}[tbh]
    \centering
    \begin{subfigure}[b]{0.49\textwidth}
        \centering
        \includegraphics[width=\textwidth]{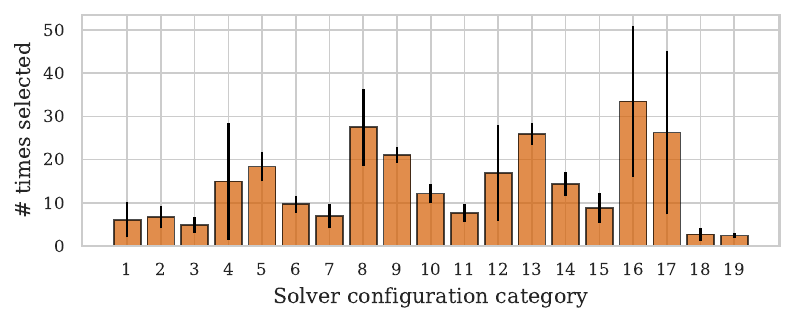}
        \caption{Heuristic}
    \end{subfigure}
    \hfill
    \begin{subfigure}[b]{0.49\textwidth}
        \centering
        \includegraphics[width=\textwidth]{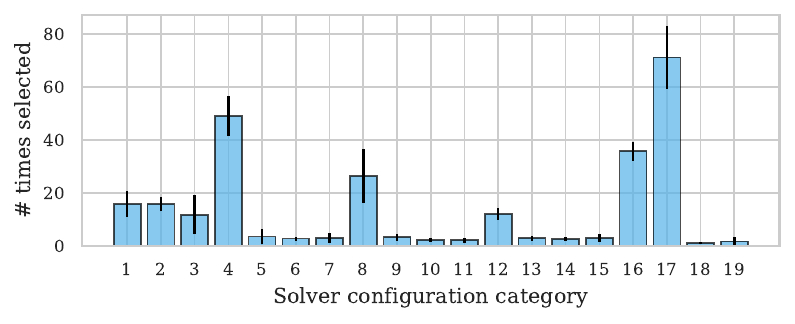}
        \caption{Gaussian Process}
    \end{subfigure}
    \caption{Histograms of the invocations of different solver configurations. Average values and standard deviations are taken w.r.t. 20 repeats of the simulation.    
    The numbering of the options corresponds to that given in Table \ref{tab:extende_solver_options}.}
    \label{fig:thermal_histogram}
\end{figure}
\begin{figure}[tbh]
    \centering
    \includegraphics[width=\textwidth]{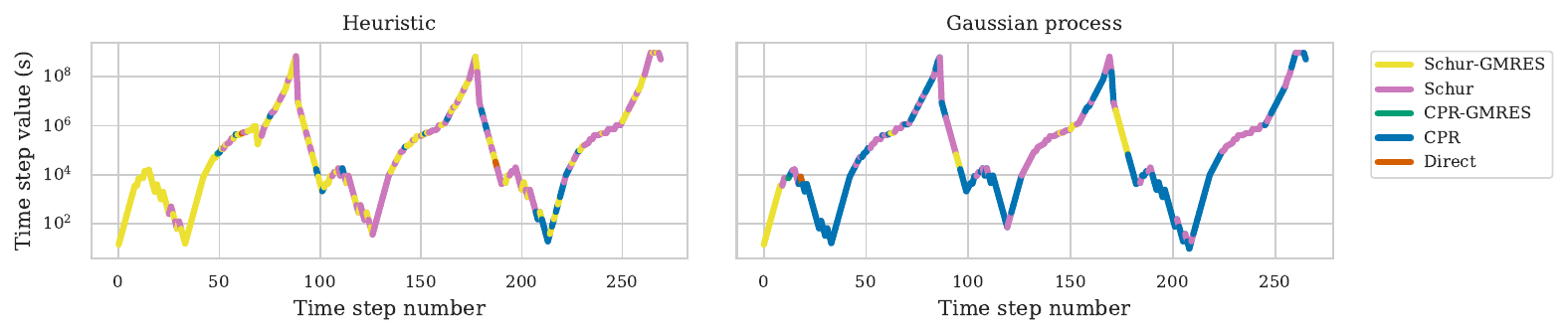}
    \caption{Time series of solver options and time step size for the heuristic and GP optimization approaches, shown for a single simulation. The solvers are grouped into five main categories defined in Section \ref{sec:experiment:extended_solver_space}.}
    \label{fig:thermal_clustered}
\end{figure}
\begin{figure}[tbh]
    \centering
    \begin{subfigure}[b]{0.49\textwidth}
        \centering
        \includegraphics[width=\textwidth]{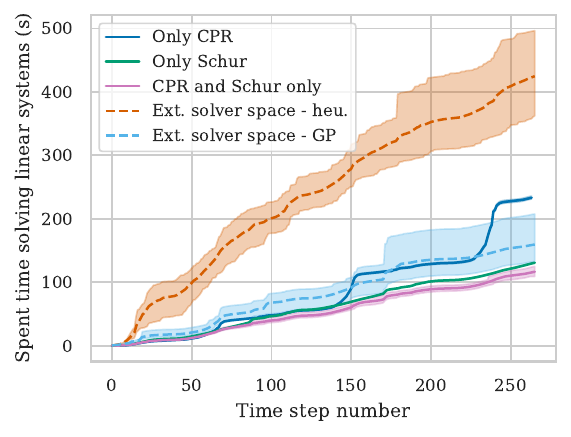}
        \caption{}
        \label{fig:many_solvers_performance}
    \end{subfigure}
    \hfill
    \begin{subfigure}[b]{0.49\textwidth}
        \centering
        \includegraphics[width=\textwidth]{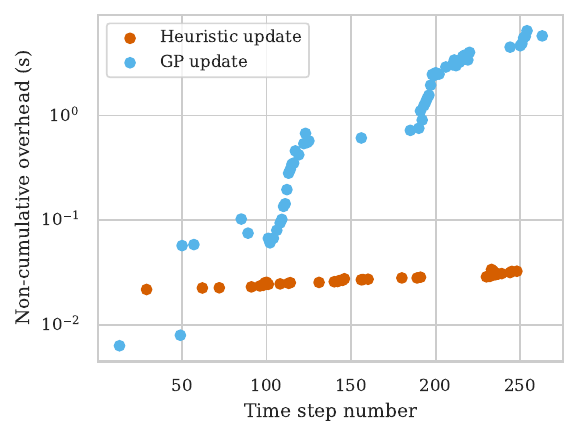}
        \caption{}
        \label{fig:overheads}
    \end{subfigure}
    \caption{(a) Cumulative performance of solver configurations picked by the solver selector in Experiment \ref{sec:experiment:extended_solver_space} for the heuristic and GP approaches (dashed lines). The data is compared with the results of the Experiment \ref{sec:numerics:categorical} (solid lines). The bold line of each color displays the mean value, while the lighter-shaded region of the same color around it indicates the range from the lowest to highest value among all the repeats of an experiment.
    (b) Overheads to train the machine learning model for the heuristic and GP approaches, shown for a single simulation and the instance of the machine learning model most commonly invoked.}
\end{figure}

Figure \ref{fig:thermal_histogram} shows a histogram of the number of invocations of the different solver configurations for the heuristic and GP optimization approaches. The data is expressed in terms of the mean and standard deviation over 20 repeated runs.
A comparison of the results for the heuristic and GP optimization methods indicates that the latter rapidly settles on a few preferred choices, while the former distributes its choices more evenly among the available options.
The heuristic method has a clear preference for the Schur complement approach, possibly reflecting that there are five times as many Schur options than CPR, thus the former will be sampled more under random exploration.
The GP optimization leads to the Schur and CPR preconditioners being called approximately equally often, and a comparison between solvers 4 and 17 (referring to Table \ref{tab:extende_solver_options}) indicates a slight preference for the implementation behind FGMRES over that of GMRES.
Furthermore, the direct solver was hardly chosen by any of the optimization methods.

To better understand the choices made, Figure \ref{fig:thermal_clustered} shows data from a single run in the same experiment, with the options grouped into five clusters: CPR and Schur with and without an inner GMRES method, and a direct solver.
The GP optimization seems to settle on a mixture of CPR and the Schur complement approach reminiscent of the optimal strategy from Section \ref{sec:numerics:categorical}, although the policy is less clear with the higher number of options.
In contrast, the heuristic method hardly considers CPR-based options at all.

The result of these policies in terms of the time spent solving linear systems is depicted in Figure \ref{fig:many_solvers_performance}.
Despite having to experiment with 19 different solvers, the average performance obtained by using GP optimization is comparable to that of CPR, and in the best case, GP approaches the pure Schur complement method. The heuristic optimization fails to find an efficient policy and makes suboptimal choices for most time steps. 
A comparison with Figure \ref{fig:thermal_performance} shows that the variations among the 20 independent runs is considerably larger than when choosing between only two options.
This holds for both the heuristic and GP optimization, and is a natural consequence of the larger number of random decisions being made during exploration.

While the GP optimization produces superior results in terms of time spent on linear solvers, this comes to the price of a significantly higher cost of the optimization problem itself.
Figure \ref{fig:overheads} depicts the time spent updating the optimization model as new data is provided after each time step.
For both the heuristic and GP methods, we show measurements for the categorical choices that were invoked most frequently, and note that these may be different between the two optimizations.
The cost of updating the heuristic model is almost constant and low, while updating the GP model is significantly more expensive, and grows with the size of the data to be fit.
We attribute the varying growth rate for the GP to the implementation of the training process, which iterative searches for the correct algorithm parameters discussed in Section \ref{sec:gp}.
The data for the other algorithmic choices is similar, although the growth in the cost of GP is less pronounced since these options have been chosen less often.
Thus, for this particular example, the total cost of solving linear systems and maintaining the optimization method is lower for the heuristic method than for the GP.
However, this observation cannot be readily transferred to other cases, as the cost of the two components scale with different factors: The time spent on the linear solver depends on the system size, coupling strength, etc., while the GP update depends on the number of data points.
The cost of the GP update can be partly alleviated, for instance, by feeding it data less frequently and in batches, or by more advanced techniques; see e.g., \cite{Vakili2021,bayesian_optimization}. 
Thus, the decision on which optimization method to use depends on the simulation setup and technology available.

\section{Conclusion}
\label{sec:conclusions}
We have presented a framework for selecting linear solvers in time-dependent multiphysics simulations. 
The framework chooses from a set of user-provided solver configurations that may include splitting schemes, preconditioners, Krylov subspace methods and numerical parameters thereof.
The selection problem is formulated as a black-box optimization problem, which we solve by Bayesian optimization with Gaussian Processes or a heuristic approach.
The key feature of our framework is that it dynamically adapts to performance data generated during the simulation and may switch solver configuration between time steps based on this information. 
This allows for quick determination of promising solver configurations without relying on pre-generated training data.
Our framework may also incorporate characteristics of the problem and the solution to dynamically adapt the linear solver.

The proposed framework was applied to two model problems, involving poromechanics and non-isothermal single-phase fluid flow in porous media, for both of which efficient linear solvers are known.
Our framework could successfully find efficient solver configurations for both problems and perform better than the pre-selected benchmarks.
The simulations also illustrated how the solver selection problem becomes harder as the number of candidate solvers grows and thus highlighted the utility in combining our framework with existing knowledge of what are efficient solver configurations.
Finally, comparison between the two optimization methods showed that the Gaussian Process gave the most accurate predictions of solver performance; however, it also had a significantly higher cost than the heuristic approach, thus the choice of optimization method is in itself non-trivial.

The task of selecting linear solvers for complex simulations is at the same time critical for performance and difficult, in particular for users of simulation software who have limited knowledge of such methods.
As our optimization approach automatizes this task, the cost and risk associated with adding options to the pool of candidate solvers is reduced, thus enabling experimentation with new and promising algorithms or implementations.
Our methodology can be inserted into existing simulation frameworks with minimal adaptation and we therefore believe it can be a useful tool to speed up numerical computations.

\section*{Acknowledgment}
This project has received funding from the VISTA program, The Norwegian Academy of Science and Letters and Equinor and from
the Norwegian Research Council, Grant 308733.
\appendix

\section{Model setup for poromechanics problem}
\label{appendix:poro}

The Mandel's problem is an established verification setup for the consolidation problem \cite{Mandel1953,Keilegavlen2021}. We solve the equation system \eqref{eq:poromechanics} in a rectangular 2D domain. 
The domain is covered by two planes from the north and south boundaries. Initially, a normal load is applied on these sides. The east and west boundaries are stress-free and allow fluid to drain. The vertical displacement at the top of the domain is time-dependent and given by the exact solution. The gravity term is neglected. The symmetry of the problem allows us to focus on the positive quarter of the domain, see \cite{Mandel1953}.

We run the Mandel's problem with the constant time step $\Delta t = 10$ s during 100 time steps. 
The values of the material parameters for the solid defined in Section \ref{sec:poromechanics} are the following:
$\mu_r=2.475 \cdot 10^9 \text{ Pa}$, $\lambda=1.65 \cdot 10^9 \text{ Pa}$, $m^{-1}=6.0606 \cdot 10^{-11} \text{ Pa}^{-1}$, $\mathbf{K} = 9.869 \cdot \mathbf{I} \cdot 10^{-14} \text{ m}^2$ (homogeneous isotropic permeability), $b=1$. For the fluid, $\rho=10^3\text{ kg}\cdot \text{m}^{-3}$, $\mu=10^{-3} \text{ Pa} \cdot \text{s}$. 
The simulated units are in SI; pressure is scaled to MPa.

\section{Model setup for non-isothermal flow in porous media}
\label{appendix:thermal}
We solve the equation system \eqref{eq:thermal_system} in a rectangular 2D domain of $(60 \times 220)$ cells with heterogeneous porosity and permeability. The porosity and permeability fields are taken from the first slice of the SPE10 benchmark \cite{SPE10}, model 2. To make permeability isotropic, we take only the first diagonal component of the permeability tensor. The permeability field is multiplied by 1000 to achieve more hyperbolic properties of heat transport. The effect of gravity is neglected. The model uses SI units, but the pressure is scaled to MPa to improve the conditioning of the system.

The values of the material parameter for the solid are the following: $\rho_r = 2650\text{ kg}\cdot \text{m}^{-3}$, $c_r=920\text{ J}\cdot \text{kg}^{-1} \cdot\text{K}^{-1}$, $k_{Tr} = 1.7295772056 \text{ W} \cdot \text{m}^{-1} \cdot \text{K}^{-1}$. For the fluid, $c_f=2093.4 \text{ J}\cdot \text{kg}^{-1} \cdot\text{K}^{-1}$.

The initial pressure in the model is $41.369$ MPa and the initial temperature is $288.706$ K. The fluid in the model is water.
All boundary conditions forbid fluid flow and energy transfer through the boundaries. The model contains two wells: one for injection and one for production. A driving force of the evolution in the model is a cyclic injection process. The model runs for 300 simulated years. Three injection periods occur during the simulation: at $t = 0, 50$ and $100$ years. Each injection period lasts 100 days; the hot fluid is injected at a temperature of $422.039$ K, a maximum injection pressure is $68.95$ MPa, and a maximum injection flux is $10^{-3}$ m$^{3}$ s$^{-1}$. The production well is limited by a minimum pressure of $27.579$ MPa and a maximum flux of $10^{-3}$ m$^{3}$ s$^{-1}$. The rest of the time, both wells are disabled.

\begin{figure}[tbh]
    \centering
    \includegraphics[width=0.5\textwidth]{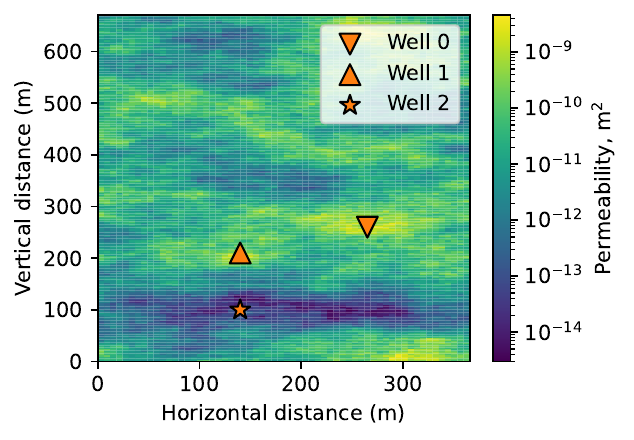}
    \caption{Well locations and permeability of the non-isothermal fluid model in \ref{sec:thermal}.}
    \label{fig:well_locations_permeability}
\end{figure}

The setup was significantly motivated by the experiment from \cite{Roy2019}, where the authors consider the Schur complement approximation preconditioner and the CPR preconditioner, described in \ref{sec:thermal_solvers}. One of the results shows that the Schur complement approximation preconditioner performs better for the original SPE10 permeability, exhibiting more parabolic behavior. On the other hand, the CPR preconditioner performs better for the permeability multiplied by 1000, exhibiting more hyperbolic behavior. We model cyclic injection to witness both solution regimes within one simulation to illustrate the case when dynamic switching of solvers is beneficial. All the constitutive laws and the material parameters are the same as in the numerical experiments in \cite{Roy2019}.

\subsubsection*{Sources}
We consider the source/sink terms that represent injection and production wells. The well locations are shown in Figure \ref{fig:well_locations_permeability}. The primarily used simulated model utilizes well 0 for injection and well 1 for production. The alternative simulated model in Section \ref{sec:experiment:knowledge_transfer} utilizes well 1 for injection and well 2 for production.
A simple way to model these is by using point sources/sinks:
\begin{equation}
    f(\mathbf{x}) = \sum_i q^i_\text{inj}(p, T) \delta(\mathbf{x - \mathbf{x^i_\text{inj}}}) \rho(p, T_\text{inj}) - \sum_j q^j_\text{prod}(p, T) \delta(\mathbf{x} - \mathbf{x}^j_\text{prod}) \rho(p, T),
\end{equation}
\begin{equation}
    f_T(\mathbf{x}) = \sum_i q^i_\text{inj}(p, T) \delta(\mathbf{x - \mathbf{x^i_\text{inj}}}) \rho(p, T_\text{inj}) cv T_\text{inj} - \sum_j q^j_\text{prod}(p, T) \delta(\mathbf{x} - \mathbf{x}^j_\text{prod}) \rho(p, T) c_v T,
\end{equation}
where $\mathbf{x_\text{inj}}$ and $\mathbf{x}_\text{prod}$ represent the location of injection and production wells, respectively. The Dirac delta $\delta(\mathbf{x})$ is used to denote the source / sink locations, $q^i_\text{inj}$ and $q_\text{prod}^j$ are injection and production volumetric rates of the wells, respectively. We consider the Peaceman well model \cite{Peaceman1978,Chen2009}. In our case of isotropic media, the rates are given by:
\begin{equation}
    q=\dfrac{2 \pi h_w K}{\mu \ln(r_e / r_w)}(p_{bh} - p),
\end{equation}
where $h_w$ is the height of well opening, $K$ is the isotropic scalar permeability, $p_{bh}$ is the bottom-hole pressure, $r_w$ is the well radius, and $r_e$ is the equivalent radius which is calculated by using the formula:
\begin{equation}
    r_e = 0.14 \sqrt{D_x^2 + D_y^2},
\end{equation}
where $D_x$ and $D_y$ are the horizontal lengths of the grid cell. We choose $h_w=5$  meters and $r_w=0.1$ meters.

\subsubsection*{Density}
The following formula is used for density:
\begin{equation}
    \rho(p, T) = \rho_0 e^{c\left(p - p_0\right)}e^{\beta\left( T - T_0 \right)},
\end{equation}
where $\rho_0$, $p_0$ and $T_0$ are reference values, $c$ is a compressibility coefficient and $\beta$ is a thermal expansion coefficient. The values of the parameters are representative of those used in reservoir simulations: $p_0=1.01325$ bar, $T_0 = 228.7056$ K, $c = 5.5 \times 10^{-5}$ bar$^{-1}$, $\beta=2.5 \times 10^{-4}$ K$^{-1}$ and $\rho_0=999$ kg m$^{-3}$.

\subsubsection*{Viscosity}
Viscosity is determined by the following correlation \cite{Bennison1998Viscosity}:
\begin{equation}
    \mu(T_F) = 10^{A_1 \gamma_{API} + A_2} T_F ^{A_3 \gamma_{API} + A_4},
\end{equation}
which takes temperature $T_F$ in 
$^{\circ}F$ and returns viscosity in cp (0.001 kg m$^{-1}$ s$^{-1}$). The values of dimensionless parameters are given:
$A_1=-0.8021$, $A_2=23.8765$, $A_3=0.31458$, $A_4=-9.21592$. The value of $\gamma_{API}$ is 10 (for water).

\bibliographystyle{unsrtnat}
\bibliography{references}

\end{document}